\documentclass{article}
\usepackage{amsmath}
\usepackage{amssymb}
\usepackage{amsthm}
\newcommand{\N}{\mathbb{N}}
\newcommand{\Z}{\mathbb{Z}}
\usepackage{enumerate}
\usepackage{hyperref}
\newtheorem{thm}{Theorem}
\newtheorem{prop}{Proposition}
\hypersetup{colorlinks=false,linkcolor=,urlcolor=black}

\title{A note on minimal additive complements of integers}
\author{Andrew Kwon\footnote{\textit{E-mail address:} \url{akwon@andrew.cmu.edu}}}
\date{}
\begin{document}
\maketitle

\abstract{Let $C, W \subseteq \mathbb{Z}$. If $C + W = \Z$, then the set $C$ is called an additive complement to $W$ in $\Z$. If no proper subset of $C$ is an additive complement to $W$, then $C$ is called a minimal additive complement. We provide a partial answer to a question posed by Kiss, S\'{a}ndor, and Yang regarding the minimal additive complement of sets of the form $W = (n \N + A) \cup F \cup G$, where $|F|<\infty, (F \mod{n}) \subseteq (A \mod{n})$ and $(G \mod{n}) \cap (A \mod{n}) = \emptyset$. We also introduce the dual problem of characterizing sets that arise as the minimal additive complements of some set of integers, proving the analog of Nathanson's initial result on existence of minimal additive complements.}\\


\section{Introduction}
Investigation of properties of sumsets of finite sets of positive integers has a long history in additive combinatorics, with primary focus being on quantitative properties such as the size of $A+A$ compared to that of $A$. However, recent focus has been placed on understanding sumsets of potentially infinite sets of integers. Nathanson \cite{Na} considered sets $W, C \subseteq \Z$ such that $W + C = \Z$. More specifically, for $W \subseteq \Z$, we say that $C \subseteq \Z$ is an \textit{additive complement} (with respect to $W$) if $C + W = \Z$ and $C$ is a \textit{minimal additive complement} (with respect to $W$) if no proper subset of $C$ is also an additive complement with respect to $W$. Nathanson posed the question of which sets $W$ have minimal additive complements. The most basic result regarding the existence of minimal additive complements is the following.

\begin{prop}\label{proposition:Nisbad} $\N$ does not have a minimal additive complement.
\end{prop}

\textit{Proof.} Consider any additive complement of $\N$, say $C$. Note that $C$ is an additive complement to $\N$ if and only if $C$ contains infinitely many negative elements. It follows that $C \setminus \{ c \}$ for any $c \in C$ will still be an additive complement to $\N$, and so $C$ is not minimal. \qed \\

Of course, the above argument can be applied to any infinite arithmetic progression that is bounded below. There are also several non-trivial classes of sets for which the existence (or nonexistence) of a minimal additive complement is known. 

\begin{thm}[Nathanson \cite{Na}] If $|W| < \infty$, then every additive complement of $W$ contains a minimal additive complement.
\end{thm}
\begin{thm}[Chen, Yang \cite{ChYa}]\label{theorem:bi-infinite} If $\inf W = -\infty$ and $\sup W = \infty$, then there exists a minimal additive complement to $W$.
\end{thm}
\begin{thm}[Chen, Yang \cite{ChYa}]\label{theorem:manycases} Let $W = \{1 = w_{1} < w_{2} < \cdots\}$ be a set of integers and 
	\[
		\overline{W} = \N \setminus W = \{\overline{w}_{1} < \overline{w}_{2} < \cdots\}.
	\]
	\begin{enumerate}[(a)]
		\item If $\limsup_{i \to \infty} (w_{i+1} - w_{i}) = \infty$, then $W$ has a minimal additive complement. 
		\item If  $\lim_{i\to\infty} (\overline{w}_{i+1} - \overline{w}_{i}) = \infty$, then $W$ does \emph{not} have a minimal additive complement.
	\end{enumerate}
\end{thm}
Theorems 1, 2, 3(a) demonstrate that types of sets that are qualitatively very different from $\N$ have minimal additive complements, while Theorem 3(b) indicates that removing an extremely sparse set from $\N$ can still yield a set that does not have a minimal additive complement. This is a part of a heuristic regarding the existence (or not) of minimal additive complements that we will further develop through this paper. With these results established, we are essentially left with considering sets of integers that are bounded below and have bounded gaps between consecutive elements.\\

In \cite{KSY}, Kiss, S\'{a}ndor and Yang consider a nice collection of such sets, namely those $W\subseteq \Z$ that are bounded below and for which there exists a positive integer $n$ such that $w + n \in W$ for sufficiently large $w \in W$, which they call \textit{eventually periodic sets}. Such sets, after some point, consist of a union of arithmetic progressions modulo $n$; the elements prior to that can be distinguished based on whether or not they are in the same residue classes modulo $n$ as the eventual arithmetic progressions. More explicitly, up to translation, all eventually periodic sets are of the form
\begin{equation}\label{eq:1}
	W = (n \N + A) \cup F \cup G,
\end{equation}
where $A \neq \emptyset$, $(F \mod n) \subseteq (A \mod n)$, $(G \mod n) \cap (A \mod n) = \emptyset$, and $|F|, |G| < \infty$. Kiss, S\'{a}ndor and Yang proved the following necessary and sufficient conditions for having a minimal additive complement for certain sets $W$ of the above form.\\

\begin{thm}[Kiss, S\'{a}ndor, Yang \cite{KSY}]\label{theorem:neccond} Let $W = (n\N + A) \cup F \cup G$ be as in \eqref{eq:1}. If $W$ has a minimal additive complement, then there exists $S \subseteq \Z/n\Z$ such that the following two conditions hold:
\begin{enumerate}[(i)]
	\item $S + (W \mod{n}) = \Z / n\Z$;
	\item For any $s \in S$, there exists $g \in G$ such that $s + (g \mod{n}) \not\in S + (A \mod{n})$.
\end{enumerate}
\end{thm}
We can interpret this result as saying that the set $G$ plays an important role in determining the existence of a minimal additive complement. In particular, if $G$ is empty then $W$ does not have a minimal additive complement. In Theorem \ref{theorem:neccond} there is no requirement for $G$ to be finite.\\

\begin{thm}[Kiss, S\'{a}ndor, Yang \cite{KSY}]\label{theorem:suffcond} Let $W = (n \N + A) \cup F \cup G$ be as in \eqref{eq:1}. Suppose there exists $S \subseteq \Z / n\Z$ such that the following two conditions hold:
\begin{enumerate}[(i)]
	\item $S + (W \mod{n}) = \Z / n\Z$;
	\item For any $s \in S$, there exists $g \in G$ such that $s + (g \mod{n}) \not\in (S \setminus \{s\}) + (W \mod{n})$.
\end{enumerate}
Then there exists a minimal additive complement to $W$.
\end{thm}

Thus, if $G$ is finite this result roughly allows us to reduce the existence of a minimal additive complement of $W$ in $\Z$ to the existence of a minimal additive complement to $W \mod{n}$ in $\Z / n\Z$. In \cite{KSY}, Kiss, S\'{a}ndor and Yang pose the question of whether Theorem \ref{theorem:suffcond} is still true when $G$ is infinite. We observe that it is certainly not true in general, for if $G$ is of the form $n \N + B$, then $W = n\N + (A \cup B)$, and such a set does not have a minimal additive complement. The first thing we do in this paper is to identify a special situation where the result does hold and provide a converse that indicates that the characteristics of the set $G$ are important for understanding the overall behavior of eventually periodic sets.\\

In the latter part of the paper, we also begin a secondary investigation into minimal additive complements that is dual in nature to the original problem posed by Nathanson: what kinds of sets of integers $C$ can arise as the minimal additive complement of some set $W$? We will then connect this line of thinking with our heuristic that minimal additive complements should arise from some sort of arithmetic irregularity by considering sets that are minimal additive complements to themselves.

\section{Statement and Discussion of Results}

In our partial answer to the question posed by Kiss, S\'{a}ndor, and Yang, we provide the following result.\\

\begin{thm}\label{theorem:infinitecase} Let $W = (n \N + A) \cup F \cup G$ be as in \eqref{eq:1}, except $G$ may be infinite and $G \pmod{n}$ consists of a single residue class. Suppose there exists $S \subseteq \Z / n\Z$ satsifying the conditions of Theorem \ref{theorem:suffcond}. If $G$ has a minimal additive complement, then so does $W$.
\end{thm}

The previous theorem indicates that, at least in the case that $G$ resides in one residue class modulo $n$, $W$ inherits the property of having a minimal additive complement from the set $G$ as long as there are no modular obstructions (which is the content of the conditions of Theorem \ref{theorem:suffcond}). We also have the following converse, which demonstrates that this is the only way $W$ will have a minimal additive complement in this case.

\begin{thm}\label{theorem:infiniteconverse} Let $W = (n \N + A) \cup F \cup G$ be as in \eqref{eq:1}, except $G$ may be infinite and $G \pmod{n}$ consists of a single residue class. If $W$ has a minimal additive complement, then so does $G$.
\end{thm}

These results together indicate that in the case of sets of this special form, one can directly identify a ``source'' of irregularity that permits a minimal additive complement to exist, namely the subset $G$. As an added consequence of these results, we can easily deduce the following theorem from \cite{KSY} without the explicit construction that was used.

\begin{thm}[Kiss, S\'{a}ndor, Yang \cite{KSY}]\label{theorem:specialcase} There exists an infinite, not eventually periodic set $W \subseteq \N$ such that $w_{i+1} - w_{i} \in \{1, 2\}$ for all $i$ and there exists a minimal complement to $W$.
\end{thm}

This follows from applying Theorem~\ref{theorem:infinitecase} for $n=2$, $A$ consisting of one element, and $G$ being a subset of one residue class that is bounded below and has arbitrarily large gaps between consecutive elements (which suffices for $G$ to have a minimal additive complement).\\

Now, all investigations regarding minimal additive complements have so far been assymetric in the sense that only the set $W$ has been studied and not the set $C$. We prove the following initial results about which sets $C$ can be minimal additive complements of some set $W$.

\begin{thm}\label{theorem:dualmacfinite} If $C \subseteq \Z$ is a finite set, then there exists a set $W \subseteq \Z$ such that $C$ is a minimal additive complement to $W$.
\end{thm}

This can roughly be thought of as the dual to Nathanson's original result from \cite{Na}. We also have the following result which is almost as obvious as Proposition \ref{proposition:Nisbad} from the beginning of this paper. 

\begin{prop}\label{prop:dualmacAP} If $C = \N$, then $C$ is not a minimal additive complement to any set $W \subseteq \Z$.
\end{prop}

Although prematurely, we adopt a similar philosophy to before, namely that a set $C$ should be qualitatively different from a infinite arithmetic progression in order to arise as the minimal additive complement of some other set $W$. Doing so allows us to have a consistent interpretation of the following interesting result at the intersection of understanding those sets that have minimal additive complements and those sets that arise as minimal additive complements.\\

\begin{thm}\label{theorem:selfmac} A set $W$ is a minimal additive complement to itself if and only if $W + W = \Z$ and $W$ contains no 3-term arithmetic progressions.
\end{thm}

Here we can begin to see connections between this theory of sumsets and the very broad combinatorial problem of avoiding arithmetic progressions. We propose to interpret this result in the following way: in order for $W$ to have a minimal additive complement, our heuristic suggests that $W$ should be quite different from an infinite arithmetic progression. However, one can imagine that the condition of arising as a minimal additive complement to some set imposes some other qualitative differences between $W$ and an infinite arithmetic progression. When both of these conditions are enforced simultaneously in this symmetric way, then we arrive at the above situation where $W$ cannot share any resemblance with an arithmetic progression at all. The above result is easy to verify, but for good measure we show it is not vacuous with the following.

\begin{prop}\label{prop:construction} There exists a set $W$ such that $W + W = \Z$ and $W$ contains no 3-term arithmetic progressions.
\end{prop}

We note that general results about sets avoiding 3-term arithmetic progressions also require that such a set $W$ must have arbitrarily large gaps between consecutive elements, which already ensures the existence of a minimal additive complement. 

\section{Proofs of Results}
We first present the proofs of our results related to understanding which sets have minimal additive complements, Theorem~\ref{theorem:infinitecase} and Theorem~\ref{theorem:infiniteconverse}.\\

\noindent \textit{Proof of Theorem~\ref{theorem:infinitecase}.} Note that if $\inf G = - \infty$, then by Theorem \ref{theorem:bi-infinite}, $W$ has a minimal additive complement. Thus, we may assume that $\inf G > - \infty$, and we further restrict our attention to the case where $G \subseteq n \N$. Let $D$ be a minimal additive complement to $G$, where $D$ must have infinitely many negative elements in each residue class of $\Z / n\Z$, and in particular every residue class in $S$. Let $\tilde{S} = \{ z \in \Z \, : \, z \pmod{n} \in S \}$ and consider $D' = D \cap \tilde{S}$. We claim that $D'$ is a minimal additive complement to $W$.\\

\par First, we show that $D'$ is a complement to $W$. In particular, we claim that
\[
	D' + (n\N + A) = \tilde{S} + (n \Z + A), D' + G = \tilde{S},
\]
and $\tilde{S} \cup (\tilde{S} + (n\Z+A)) = \Z$ so that $D' + W = (D' + (n\N + A)) \cup (D' + G)  = \Z$. The latter equation follows from the observation that if 
\[
	D_{i} = \{ d \in D \, : \, d \equiv i \pmod{n}\}
\]
for $1 \leq i \leq n$, then the sumset $D_{i} + G$ will be precisely $n \Z + i$ since 
\[
	\bigcup_{i=1}^{n} (D_{i} + G) = D + G = \Z
\]
and $D_{i} + G \subseteq n \Z + i$. Thus, $D' + G$ will consist of all integers whose residues modulo $n$ are contained in $S$, which is precisely $\tilde{S}$. For the former equation, we claim that every integer whose residue modulo $n$ is in $S + (A \mod{n})$ will be in the sumset $D' + (n\N + A)$. Indeed, for any residue $x \in S + (A \mod{n})$ there is some integer $m$ with $m \mod{n} = x$ so that there exist $d \in D', a \in A, \ell \in \N$ such that 
\[
	m = d + (\ell n + a).
\]
Furthermore, $m +k n = d' + ((\ell + k)n + a) \in D' + (n \N + A)$ for all $k \geq 0$. To represent $m - k n, k \leq 0$ as an element of $D' + (n \N + A)$, we can write 
\[
	m - k n = d' + (\ell' n + a)
\]
for some sufficiently negative $d' \in D'$ with the same residue modulo $n$ as $d$ and some appropriately chosen $\ell' \in \N$. The choice of $d'$ is possible because $D'$ contains infinitely many negative elements in each of its residue classes in $D' \mod{n}$.\\

Finally, since $S + (W \mod{n}) = (S + (A \mod{n})) \cup C = \Z / n\Z$, we find that 
\[
	(\tilde{S} + (n \Z + A)) \cup \tilde{S} = \Z,
\]
and so $D' + W = \Z$.\\

\par To conclude that $D'$ is minimal, it suffices to show that $S + (A \mod{n})$ and $S + (G \mod{n}) = S$ are disjoint in $\Z/n\Z$. Indeed, if this is the case then for any $d \in D'$ there exists a $g \in G$ such that $d + g$ is uniquely represented in $D' + G$ since $d \in D$ and $D$ is a minimal additive complement. Also, $d + g$ is not in any of the residue classes represented by $D' + (n \N + A)$, so it is uniquely represented in $D' + W$ as well. To finish, we note that the conditions on $S$ require that for all $c \in C$, 
\[
	s \not\in (S \setminus \{s\}) + (W \mod n),
\]
and so in particular $s \not\in (S \setminus \{s\}) + (A \mod n)$. Furthermore, $s \not\in s + (A \mod n)$ since we presume that $G \mod n = \{0\}$ and thus $0 \not\in A \mod n$. Therefore, $s \not\in S + (A \mod n)$. We conclude that $W$ has $D'$ as a minimal additive complement when $G \subseteq n \N$.\\

To conclude the general case, we observe that the above argument is unaffected by translation of the set $W$, and so we can translate the set $W$ so that $G \subseteq n \N$, at which point we are done. \qed  \\

\noindent \textit{Proof of Theorem~\ref{theorem:infiniteconverse}.} As before, we assume for simplicity that $G \subseteq n \N$. Suppose $C$ is a minimal additive complement to $W$, and consider the partition $C = C_{1} \cup C_{2}$, where $C_{1}$ is the restriction of the set $C$ to the residue classes modulo $n$ in which $C$ has infinitely many negative elements, and $C_{2} = C \setminus C_{1}$. Then, 
\begin{align*}
	W + C &= (n\N + A + C) \cup (F + C) \cup (G + C)\\
	&= (A + C_{1} + n \Z) \cup (A + C_{2} + n \N) \cup (F + C_{2} + n \N) \cup (G + C),
\end{align*}
where $A + C_{1} + n \N = A + C_{1} + n \Z$ since $C_{1}$ contains infinitely many negative elements modulo $n$. Furthermore, every element of $F$ is smaller than the element of $A$ that is in the same residue class modulo $n$, and so $A + C_{2} + n \N \subseteq F + C_{2} + n \N$. We are thus left with 
\[
	\Z = W + C = (A + C_{1} + n\Z) \cup (F + C_{2} + n \N) \cup (G + C).
\]
One can imagine that $A + C_{1} + n \Z$ is a union of some number of dilated and translated ``copies'' of $\Z$, completely filling some residues classes modulo $n$, while $F + C_{2} + n \N$ is some number of dilated and translated ``copies'' of $\N$, filling half of some residue classes modulo $n$. The remainder of the argument is based on analyzing how $G + C_{1}$ interacts with these copies of $\Z$ and $\N$.\\

Say an integer $z \in \Z$ is \emph{dependent} on $c \in C$ if $z \not\in (C \setminus \{c\}) + W$. That is, $z$ is one of the integers that $c$ plays a role in uniquely representing in the sumset $C + W$. In particular, we necessarily have $z \in c + W$. Now, evidently no integer in $A + C_{1} + n\Z$ will be dependent on any element of $C$, since every integer in that set can be written as an infinite number of pairs $a + c + nm$, where $a + nm \in A + n \N$ and $c \in C_{1}$. Thus, all of the integers dependent on some $c \in C_{1}$ must be contained in $G + C_{1}$. From this we also conclude easily that $G + C \not\subseteq A + C_{1} + n \Z$. There are two cases to consider.\\

Suppose that $(A + C_{1} \cup F + C_{2}) \pmod{n} \neq \Z / n \Z$. That is, some residue class $i + n \Z$ is disjoint from $A + C_{1} + n\Z$ and $F + C_{2} + n \N$. Then, $i + n \Z$ must be contained in $G + C_{1}$, since $G + C_{2}$ does not contain any arbitrarily negative elements (and $G + C_{1}, G + C_{2}$ occupy different residue classes modulo $n$). Letting $C_{1,i} = C_{1} \cap (i + n \Z)$ and using the assumption that $G \subseteq n \N$, it is not hard to see that $G + C_{1,i}$ is the only subset of $G + C_{1}$ that meets $i + n\Z$, while $G + C_{1}$ contains $i + n\Z$, and so we must have
\[
	G + C_{1,i} = i + n\Z.
\]
Furthermore, every $c \in C_{1,i}$ has some $z \in G + C_{1}$ that is dependent on $c$, but of course such a $z$ is more specifically contained in $G + C_{1,i}$. Thus, no proper subset $C'$ of $C_{1,i}$ can satisfy $G + C' = i + n\Z$. By taking 
\[
	D = \bigcup_{j=0}^{n-1} (C_{1,i} + j),
\]
we find that $G + D = \Z$ and $D$ is a minimal additive complement because every element of $C_{1,i} + j$ has some $z \in i + j + n \Z$ that is dependent upon it. Thus, $G$ has a minimal additive complement in this case.\\

The other case to consider is that $(A + C_{1} \cup F + C_{2}) \pmod{n} = \Z / n \Z$. Evidently it is impossible for $A + C_{1} \pmod{n} = \Z / n \Z$, since then $(n \N + A) + C_{1} = Z$, and so $C$ cannot be a minimal additive complement to $W$. Thus, the restriction of $G + C_{1}$ to some residue class $i + n \Z$ will meet $F + C_{2} + n \N$ but not $A + C_{1} + n \Z$. This restriction will be precisely $G + C_{1,i}$, where $C_{1,i}$ is defined as above. Furthermore, all of the integers $z \in i + n \Z$ that are dependent on some $c \in C_{1,i}$ must be less than the minimum element of $(F + C_{2} + n \N) \cap (i + n \Z)$, say $m$. This follows from the observation that all integers greater than $m$ that are in $G + C_{1,i}$ will already be contained in $F + C_{2} + n \N$, and therefore are not dependent on any elements of $C_{1,i}$. Considerations of residue classes similar to before guarantee that $G + C_{1,i}$ will contain all sufficiently negative integers, specifically all that are less than $m$.\\

Thus, we are in a situation where $G + C_{1,i}$ contains all integers in $i + n \Z$ less than some integer $m$, and all integers $z$ dependent on some $c \in C_{1,i}$ must also be less than $m$. Now, if $G + C_{1,i} = i + n \Z$, then we are done by the above argument; otherwise, consider the smallest integer $z > m$ in $i + n \Z$ that is not in $G + C_{1,i}$. Letting $g = \min G$, instead of the set $C_{1,i}$ we can instead consider the set $C_{1,i} \cup \{ z - g\}$. Then, $z$ will be dependent on $z-g$, and adding the element $z-g$ will not affect any of the integers dependent on elements of $C_{1,i}$, and so we can iterate this process to inductively define a set $C_{i}$ such that $G + C_{i} = i + n \Z$ where every element of $c \in C_{i}$ has some integer $z \in i + n \Z$ such that $z$ is dependent on $c$. Then we are in the same situation as in the previous case, and so we again can conclude that $G$ has a minimal additive complement. \qed\\

Now we present the proofs of our results related to the dual problem of understanding which sets can arise as minimal additive complements, Theorem~\ref{theorem:dualmacfinite}, Proposition~\ref{prop:dualmacAP}, Theorem~\ref{theorem:selfmac}, and Proposition~\ref{prop:construction}.\\

\noindent \textit{Proof of Theorem~\ref{theorem:dualmacfinite}.} Let $C = \{c_{1}, \ldots, c_{n}\}$ with $c_{1} < c_{2} < \cdots < c_{n}$ and $k$ be the maximum difference between consecutive elements of $C$. Then, it is easy to see that 
\[
	C + [k+1] = \{ z \, : \, c_{1} + 1 \leq z \leq c_{n} + k+1\}.
\]
Indeed, for any $z \in [c_{1} + 1, c_{n}]$, there are $c_{i}, c_{i+1}$ such that $c_{i} \leq z < c_{i+1}$. Then, since $c_{i+1} - c_{i} < k+1$, there must be $j < k+1$ such that $c_{i} + j = z$. Therefore $[c_{1} + 1, c_{n}] \subseteq C+[k+1]$. For $z \in [c_{n} + 1, c_{n} + k+1]$, it is obvious that $z \in C + [k+1]$ and so $[c_{n} + 1, c_{n} + k + 1] \subseteq C + [k+1]$. Clearly $C + [k+1] \subseteq [c_{1} + 1, c_{n} + k + 1]$, and so we have equality, as desired.\\

Now, we will outline an iterative construction for a set $W \subseteq \Z$ such that $W + C = \Z$ with that $C$ being a minimal additive complement to $W$. Start with $W = \{z \, : \, z \leq -c_{n} - 1\}$. With this initial choice, it is not hard to see that $W + C$ is precisely the set of negative integers. We now use the following procedure:\\

\noindent For each $i \in [n]$:\\
\indent If $z_{i}$ is the least integer currently not in $W+C$, include $z_{i} - c_{i}$ and $[k + 1] + (z_{i} - c_{1}) = [z_{i} + 1 - c_{1}, z_{i} + k + 1 - c_{1}]$ in $W$. \\

In each step of the above, $z_{i}$ will certainly be in $W + C$ since $(z_{i}-c_{i}) + c_{i} = z_{i}$. On the other hand, $C + [k + 1] + (z_{i} - c_{1})$ will be the set of all integers in $[z_{i} + 1, z_{i} + c_{n} + k + 1 - c_{1}]$. We aim to show that this $z_{i}$ will only be represented once in $W + C$, by $c_{i} \in C$ and $z_{i} - c_{i} \in W$. To deduce this, it suffices to show that every subsequent element $w$ we add to $W$ will be too large for $C + w$ to even come close to $z_{i}$.\\

Indeed, immediately after iteration $i$ of the above, the next integer not in $W + C$ will be $z_{i+1} = z_{i} + c_{n} + k + 2 - c_{1}$, and so the elements we add to $W$ in iteration $i+1$ will be at least
\[
	z_{i} + c_{n} + k + 2 - c_{1} - c_{i+1}.
\]
However, the minimum element of $C + (z_{i} + c_{n} + k + 2 - c_{1} - c_{i+1})$ is 
\[
	c_{1} + (z_{i} + c_{n} + k + 2 - c_{1} - c_{i+1}) = z_{i} + c_{n} - c_{i+1} + k + 2 > z_{i},
\]
and so all of the translates of $C$ that arise from adding new elements to $W$ according to this procedure are bounded away from $z_{i}$. Therefore, after applying this procedure we will have found a set $W$ such that for each $c_{i} \in C$, there exists an integer $z_{i}$ such that $z_{i} \not\in W + (C \setminus \{c_{i}\})$. Lastly, we observe that at this point $W + C$ is precisely the set of integers at most $z_{n} + c_{n} + k + 1 - c_{1}$, and so we can simply include many large translates of $[k+1]$ to $W$ in order to get $W + C = \Z$. \qed\\

\noindent \textit{Proof of Proposition~\ref{prop:dualmacAP}.} Suppose that there is a set $W \subseteq \Z$ such that $\N$ is a minimal additive complement to $W$. That is, $W + \N = \Z$, and for every $n \in \N$ there is some $z \in \Z$ such that $z \not\in W + (\N \setminus \{n\})$. However, in order for $W + \N$ to equal $\Z$, $W$ must be unbounded from below, i.e., contain infinitely many negative elements. Therefore, for $w \in W$ that are sufficiently negative, we can find that $z - w \neq n$, which is a contradiction. Therefore $\N$ cannot be a minimal additive complement to any set. \qed\\

\noindent \textit{Proof of Theorem~\ref{theorem:selfmac}.} First, suppose $W$ contains some arithmetic progression $w-d, w, w+d$ with $d \neq 0$. We will show that $W$ cannot be its own minimal additive complement. Suppose to the contrary that it was; then, there must be some $z$ such that $z \not\in W + (W \setminus\{w\})$. However, if $z = v + w$ with $v \neq w$, then $z$ will be represented in $W + (W \setminus\{w\})$, simply as $w + v$ instead of $v + w$. Therefore, $v = w$ and $z = 2w$. However, in this case we also see the obstruction arising from arithmetic progressions, as $(w - d) + (w + d) = 2w = z$. Thus $W$ is not a minimal additive complement to itself.\\

The proof of the opposite direction is even easier. If $W$ contains no 3-term arithmetic progressions, then for each $w \in W$, $2w$ is uniquely represented in $W + W$ as $w + w$, thus the result.\qed\\

\noindent \textit{Proof of Proposition~\ref{prop:construction}.} We use a construction inspired by the classical greedy construction for the set of integers that does not contain any 3-term arithmetic progressions, namely the set of nonnegative integers with only 0s and 1s in their ternary expansions. Let $A$ denote this set, and $2A = \{2a \, : \, a \in A\}$. Furthermore, it is not hard to see that $A + A = \N$. Then, $2A + 2A = 2\N$, the set of all even nonnegative integers. Our choice of $W$ will be $W = 2A \cup (-2A - 1)$.\\

Notice that this set certainly avoids 3-term arithmetic progressions; there are no 3-term arithmetic progressions in $2A$ or $-2A - 1$, and parity considerations ensure there will not be any 3-term arithmetic progressions with terms in both $2A, -2A - 1$. Now, in considering $W + W$, we can separately analyze contributions from essentially three distinct sumsets, $2A + 2A$, $2A + (-2A - 1)$, and $(-2A - 1) + (-2A - 1)$. As mentioned above, $2A + 2A$ will yield all even nonnegative integers, and similarly $(-2A - 1) + (-2A -1)$ will yield all even negative integers. Now we direct our focus on $2A + (-2A - 1)$.\\

We first consider positive odd integers $2z + 1$. We would like to exhibit the existence of $a_{1}, a_{2} \in A$ such that $2a_{1} - 2 a_{2} - 1 = 2 z + 1$, i.e. $2z + 1 \in 2A + (-2A - 1)$. We can rearrange this equation into $a_{1} - a_{2} = z + 1$, so this simplifies to the following: given an arbitrary positive integer $z$, is there an $a \in A$ such that $a - z \in A$? We outline a means by which to construct such an $a$.\\

Let $z_{-1} = z$ and define $z_{j}$ in the following way for $j \geq 0$. If the $j^{\text{th}}$ least significant ternary digit of $z_{j-1}$ is 2, define $z_{j} = z_{j-1} + 3^{j}$. Otherwise, take $z_{j} = z_{j-1}$. It is not hard to see that this procedure will ensure that at step $j$, none of the $j$ least significant digits will be 2. Furthermore, this sequence will be eventually constant (since eventually all of the ternary digits will be 0), say with value $\tilde{z}$. We claim that $\tilde{z} \in A$ and $\tilde{z} - z \in A$. That $\tilde{z} \in A$ is clear, while $\tilde{z} - z$ consists only of lone powers of 3 that were added when necessary; no power of 3 was added twice, and so this observation is clear as well.\\

Thus, given any $z \geq 0$ we can find $a_{1}, a_{2}$ such that $a_{1} - a_{2} = z + 1$, and so all odd positive integers are in $2A + (-2A - 1)$. Now, deducing the presence of all odd negative integers is easy; given $a_{1}, a_{2}$ with $a_{1} - a_{2} = z+1$ we can reverse the roles of $a_{1}, a_{2}$ to find 
\[
	2 a_{2} - 2 a_{1} - 1 = - 2z + 1,
\]
which for $z \geq 1$ exactly parametrizes the odd negative integers.\\

We conclude that this choice of $W$ avoids 3-term arithmetic progressions and also satisfies $W + W = \Z$, as desired. \qed

\section{Concluding Remarks}
Theorems~\ref{theorem:infinitecase},\ref{theorem:infiniteconverse} provide the first results demonstrating that certain special subsets control the existence of a minimal additive complement or not. However, extending beyond the case where $G$ resides in a single residue class is difficult because it is no longer clear which subset matters, whether it could be the entirety of $G$, or some subset of it (plausibly restriction of $G$ to some residue class with important properties). There are many other natural questions stemming from this line of inquiry; in the work of Kiss, S\'andor, and Yang (as well as the work here), the set $F$ is entirely inconsequential, but it is not understood in general whether adding or removing elements from a set changes the existence of a minimal additive complement or not. A possibly more pressing issue is in discovering more sufficient conditions for a set to not have a minimal additive complement, for which there is only the one result due to Chen and Yang.\\

There are also many questions in the newer line of inquiry regarding which sets arise as minimal additive complements. The author expects results similar to Theorems~\ref{theorem:bi-infinite},\ref{theorem:manycases} should be attainable. With or without that additional information, it could be interesting to investigate what kinds of properties can be deduced about sets $W$ that have minimal additive complements and also arise as minimal additive complements. The very special case where $W$ is its own minimal additive complement reduced to a rather simple characterization that could be roughly explained by heuristics discussed in this paper, but investigating how those two properties interact in general may yield much deeper insights.

\section{Acknowledgements}
This research was conducted at the University of Minnesota Duluth REU and was supported by NSF grant 1659047. The author thanks Joe Gallian for suggesting the problem and conducting the REU, Levent Alpoge and Sean Nemetz for helpful discussions, and Ben Gunby and Colin Defant for their feedback on the manuscript.

\end{document}